\documentclass[10pt]{amsart}
\usepackage{amssymb,amscd}
\usepackage[all,cmtip]{xy}
\usepackage{graphics}
\usepackage{epic}
\usepackage{url}
\usepackage[english]{babel}
\usepackage[utf8]{inputenc}
\usepackage{multirow}

\newtheorem{theorem}{Theorem}[section]
\newtheorem{lemma}[theorem]{Lemma}
\newtheorem{proposition}[theorem]{Proposition}

\newtheorem{definition-lemma}[theorem]{Definition-Lemma}

\newtheorem{definition}[theorem]{Definition}

\newtheorem{example}[theorem]{Example}

\newtheorem{remark}[theorem]{Remark}

\def\CC{{\mathbb C}}

\def\EE{{\mathbb E}}
\def\FF{{\mathbb F}}

\def\QQ{{\mathbb Q}}

\def\ZZ{{\mathbb Z}}
\def\cO{{\mathcal O}}

\def\cE{{\mathcal E}}
\def\cD{{\mathcal D}}

\newcommand\Hom{\operatorname{Hom}}

\newcommand\Ker{\operatorname{Ker}}

\newcommand\Spec{\operatorname{Spec}}

\begin{document}

\title[Generalization of Deuring Reduction Theorem]{Generalization of Deuring Reduction Theorem}

\author{Alexey Zaytsev}

\email{alzaytsev@kantiana.ru}
\address{I. Kant Baltic Federal University, Nevskii 14, Kaliningrad, Russia}
%%%%%%%%%%%%

\maketitle

 % \today
\begin{abstract}
In this paper we generalize the Deuring theorem on a reduction of elliptic curve with complex multiplication. More precisely, for an Abelian variety $A$, arising after reduction of an Abelian variety with complex multiplication by a CM field $K$ over a  number field
at a pace of good reduction. We establish a connection between a decomposition of the first truncated Barsotti-Tate group scheme $A[p]$ and a decomposition of $p\cO_{K}$ into prime ideals. In particular, we produce these explicit relationships  for Abelian varieties of dimensions $1,\, 2$ and $3$.
\end{abstract}

\begin{section}{introduction and results}\label{introduction}

In this paper we generalize the Deuring reduction theorem on elliptic curves with complex multiplication to Abelian varieties with complex multiplication. 

We start with the  Deuring reduction theorem and its reformulation in terms of group schemes.
The classical Deuring reduction theorem (see for example \cite{Elliptic_Functions}, theorem~12, page 182) states:
\begin{theorem}\label{Deuring}
Let $\cE$ be an elliptic curve over a number field, with ${\rm End}(\cE)\cong \cD$,
where $\cD$ is an order in an imaginary quadratic field $K$. Let $\mathcal{P}$ be a place of 
$\bar{\QQ}$ over a prime number $p$, where $\cE$ has non-degenerate reduction $E$.
The curve $E$ is supersingular if and only if $p$ has one prime of $K$ above it (p is ramified or inert). The curve $E$ is ordinary if and only if $p$ splits completely in $K$.
\end{theorem}
In other words, as an abstract group the group of $p-$torsion points of an elliptic curve $E$ is isomorphic to
$$
E[p](\bar{\FF}_p)\cong\left\{
\begin{array}{l}
(0), \, \mbox{if}\quad p\cO_{K}=\mathcal{P}^2\quad \mbox{or}\quad p\cO_{K}=\mathcal{P},\\
\ZZ/p\ZZ, \quad \mbox{if}\quad p\cO_{K}=\mathcal{P} \mathcal{P}^{c},\\
\end{array}
\right.
$$
where $P^{c}$ is the complex conjugation of $P$.

Furthermore, this theorem can be reformulated in a modern language, namely in terms  of group schemes or more precisely in terms of the first truncated Barsotti-Tate group schemes,
$$
E[p]\times \bar{\FF}_p\cong\left\{
\begin{array}{l}
{\rm I}_{1,1}, \quad \mbox{if}\quad p\cO_{K}=\mathcal{P}^2\quad \mbox{or}\quad p\cO_{K}=\mathcal{P},\\
\ZZ/p\ZZ \times \mu_{p}, \quad \mbox{if}\qquad p\cO_{K}=\mathcal{P} \mathcal{P}^{c},\\
\end{array}
\right.
$$
where  ${\rm I}_{1,1}$ is an indecomposable  (but not simple) group scheme of order $p^2$ fitting into an exact sequence
$$
0 \rightarrow \alpha_p \rightarrow {\rm I}_{1\, 1} \rightarrow  \alpha_p \rightarrow  0,
$$
and 
\begin{displaymath}
\begin{array}{c}
\alpha_{p}={\rm Spec}(\bar{\mathbb{F}}_{p}[X]/(X^p)) \cong {\rm Ker}({\rm Frob}: \mathbb{G}_{a} \rightarrow \mathbb{G}_{a}),\\
\mu_{p}={\rm Spec}(\bar{\mathbb{F}}_{p}[X]/(X^p-1)) \cong {\rm Ker}({\rm Frob}: \mathbb{G}_{m} \rightarrow \mathbb{G}_{m}).\\
\end{array}
\end{displaymath}

In Section~\ref{KraftRT} it is shown how a decomposition of  ideal $p\cO_{K}$ into primes implies a decomposition of the first de Rham cohomology ${\rm H}^1_{\rm DR}(A/\bar{\FF}_p)$ and an action of the Frobenius on it. Due to the equivalence of categories of the first de Rham cohomology ${\rm H}^1_{\rm DR}(A/\bar{\FF}_p)$ and the category of Dieudonn\'e modules of   $A[p] \times \bar{\FF}_p$ and  Kraft diagrams, we are able to derive a Kraft diagram of  $A[p] \times \bar{\FF}_p$ and hence an isomorphism class of the ${\rm BT}_1$-group scheme. 

As consequence of this relationship we obtain a generalization of the Deuring theorem on Abelian varieties of arbitarary dimension on the base of a Galois group of a Galois closure $\tilde{K}$ of $K/ \QQ$, a CM-type and a generator of the decomposition group of an ideal in $\tilde{K}$ lying over a prime number $p$. We carry out this procedure for  dimensions $1,\, 2$ and $3$ relying on the paper \cite{CM-Fields} and the classification of ${\rm BT}_1$-group schemes in order to obtain the following theorem.

\begin{theorem}\label{Main_Theorem}
Let $\mathcal{A}$ be an abelian scheme over ${\rm Spec}(\cO_{L})$, where $\cO_{L}$ is the ring of integers of a number field $L$. Assume that 
\begin{itemize}
\item{$  \cO_{K} \hookrightarrow {\rm End}(\mathcal{A})$, with $\cO_{K}$ is  the full ring of integers of a CM field $K$ and $K \subset L$, } 
\item{$\mathcal{A}$ has a good reduction at prime $\mathcal{P}$ of $\cO_{L}$, denoted by $A$,}
\item{$\mathcal{A}$ is simple (or, equivalent, CM-type of $\mathcal{A}$ is primitive),}
\item{a prime $p=\mathcal{P} \cap \ZZ$ is unramified in $K$.}
\end{itemize}
Then the following relations between decomposition of $p\cO_{K}$ into primes, decomposition of $BT_{1}-$group scheme $A[p]$, $p$-rank and  $a$-number hold.

\begin{tabular}{|c|c | c | c | c |}
\hline
{\sf ideal decomposition} &{\sf $BT_{1}-$group scheme $A[p]$} &{\sf  $p$-rank} &{\sf $a$-number} \\
\hline
\multicolumn{4}{|c|}{\bf  Elliptic Schemes}\\
\hline
$PP^{c}$ & $ \ZZ/p\ZZ \times \mu_{p}$ & $1$ & $0$ \\
\hline
$P$ & ${\rm I}_{1,1}$& $0$&$1$ \\
\hline
\multicolumn{4}{|c|}{\bf Abelian Surfaces}\\
%\hline
%ideal decomposition &$BT_{1}$ group schemes & $p$-rank & $a$-number\\
\hline
$P_{1}P_{1}^{c}P_{2}P_{2}^{c}$& $(\ZZ/p\ZZ \times \mu_{p})^{2}$ & $2$&$0$\\
\hline
\multirow{2}{*}{$P P^{c}$ }&   $(\ZZ/p\ZZ \times \mu_{p})^{2}$ & $2$&$0$\\ \cline{2-4}
&   ${\rm I}_{1,1}^{2}$ & $0$&$2$\\
\hline
$P_{1}P_{1}^{c}P_{2}$ & $\ZZ/p\ZZ \times \mu_{p} \times {\rm I}_{1,1}$ & $1$&$1$\\
\hline
$P_{1} P_{2}$ & ${\rm I}_{1,1}^{2}$ &$0$ &$2$\\
\hline
$P$ & ${\rm I}_{2,1}$ & $0$&$1$\\
\hline
\multicolumn{4}{|c|}{\bf Abelian Threefolds}\\
%\hline
%ideal decomposition &$BT_{1}$ group schemes & $p$-rank & $a$-number\\
\hline
$P_{1} P_{1}^{c}P_{2} P_{2}^{c}P_{3} P_{3}^{c}$& $(\ZZ/p\ZZ)^{3} \times (\mu_{p})^{3}$ & $3$ & $0$\\
\hline
\multirow{2}{*}{$P_{1} P_{1}^{c}P_{2} P_{2}^{c}$} &$(\ZZ/p\ZZ)^{3} \times (\mu_{p})^3 $& $3$&$0$\\ \cline{2-4}
& $\ZZ/p\ZZ \times \mu_{p} \times {\rm I}_{1,1}^{2}$ & $1$ &$2$\\
\hline
\multirow{2}{*}{$PP^{c}$}& $(\ZZ/p\ZZ)^{3} \times (\mu_{p})^3 $ &$3$&$0$\\\cline{2-4}
&${\rm I}_{3,2}$ &$0$&$2$\\
\hline
$P_{1}P_{1}^{c}P_{2}P_{2}^{c}P_{3}$ &$(\ZZ/p\ZZ\times \mu_{p})^{2}\times {\rm I}_{1,1}$ &$2$&$1$\\
\hline
\multirow{2}{*}{$P_{1}P_{1}^{c}P_{2}$}& $(\ZZ/p\ZZ \times \mu_{p})^{2} \times {\rm I}_{1,1}$ & $2$&$1$\\ \cline{2-4}
& ${\rm I}_{1,1}^{3}$ &$0$&$3$\\
\hline
\multirow{2}{*}{$P_{1}P_{1}^{c}P_{2}P_{3}$} &$(\ZZ/p)^{3} \times (\mu_{p})^{3}$&$3$&$0$\\ \cline{2-4}
& $\ZZ/p\ZZ \times \mu_{p} \times {\rm I}_{1,1}^{2}$ &$1$&$2$\\
\hline
$P_{1}P_{2}P_{3}$ & ${\rm I}_{1,1}^{3}$ &$0$&$3$\\
\hline
$P_{1}P_{2}$ & ${\rm I}_{2,1} \times {\rm I}_{1,1}$ & $0$&$2$\\
\hline
\multirow{2}{*}{$P$} & ${\rm I}_{3,1}$ &$0$ &$1$\\ \cline{2-4}
& ${\rm I}_{1,1}^{3}$& $0$&$3$\\
\hline
\end{tabular}
\end{theorem}

%%%%%%%%%%%%%%%%%%%%%%%
\begin{remark}
If $A$ is an Abelian variety  and
$R$ is an order of a CM field $K$ with
$[K:\QQ]=2\,{\rm dim}(A)$ such that
 there is an invective homomorphism $R\hookrightarrow {\rm End}(A)$, 
 then we say that $A$ has {\bf a complex multiplication by a ring $R$} 
 or    {\bf a complex multiplication by a CM field $K$}.
 \end{remark}

At the end of this paper we attach  a new parameter, $b-$number, to each isogeny class of ordinary Abelian varieties. In the same section a notion of density is introduced and calculated for dimension two over some finite fields.
\begin{subsection}*{Acknowledgment}
I would like to thank professor Frans Oort for and professor Brian Conrad for  their kind support and attention to my work. 
I also would like to express my deep gratitude to professor Torsten Ekedahl for explaining me his paper which is a basis of my work. 
\end{subsection}

\end{section}
%------------------------------------------------------------------------------------
\begin{section}{category of finite commutative group schemes over perfect field}
In this section we give a brief description of the category of finite commutative group schemes over a perfect field and its basic properties which are needed.  More detail information can be found in \cite{group_schemes_Pink}, \cite{group_schemes_Oort},
\cite{Oort_simple}, \cite{Moonen_group_schemes}, \cite{Demazure}  and \cite{Kraft}.

\begin{subsection}{Decomposition of the category}
Let $k$ be a perfect field of characteristic $p>0.$ 
Denote by $C=C_{k}$ the category of finite commutative group schemes over $\Spec(k)$. 
By the definition a finite scheme $G$ over $\Spec(k)$ is affine. Therefore the category $C$ is equivalent to the category of commutative finitely generated $k$-bialgebras (which are automatically flat).

\begin{proposition}
The connected component $G^{0}$ of a zero section in a finite commutative group scheme $G$ is a closed finite commutative subgroup scheme and the quotient $G/G^{0}$ is \'etale.
\end{proposition}
 
 Now we can introduce the following classes of finite commutative group schemes.
\begin{definition}
A finite commutative group scheme $G$ over $\Spec(k)$ is called 
\begin{itemize}
\item{{\bf \'etale} if the structure morphism $G \to \Spec(k)$ is \'etale} 
\item{{\bf local} if $G$ is connected} 
\end{itemize}
\end{definition}

This definition can be reformulated in terms of Hopf algebras, namely
a commutative finite group scheme $G=\Spec(A)$ is 
\begin{itemize}
\item{{\bf \'etale} if $A$ is a separable algebra over $k$ (in particular, $A$ is reduced),} 
\item{{\bf local} if $A$ is  a local ring.} 
\end{itemize}

We will write $C_{loc}$ for the full subcategory of $C$ consisting of all $G \in C$ which are local and $C_{et}$ for the full subcategory of $C$ consisting of all $G \in C$ which are \'etale.

The important observation is the following.  
\begin{lemma}
Let $G \in G_{loc}$ and $H \in C_{et}$. Then
\[
{\rm Hom}_{C}(G,H)={\rm Hom}_{C}(H,G)=(0).
\]
\end{lemma}

Each object $G \in G$ can written in a unique way in the form
$$
G=G_{et}\times G_{loc},
$$
where $G_{et} \in C_{et}$, $G_{loc} \in C_{loc}$ and 
\[
{\rm Hom}_{C}(G,H)={\rm Hom}_{C}(G_{et},H_{et}) \times {\rm Hom}_{C}(G_{loc},H_{loc}).
\]

The same decomposition holds for the linear dual $A^{*}$. Therefore the category $C$
 splits into four categories:
\begin{itemize}
\item{$C_{loc, loc}$ the category of all  $G \in C_{loc}$ with $G^D \in C_{loc}$,}
\item{$C_{loc, et}$ the category of all  $G \in C_{loc}$ with $G^D \in C_{et}$,}
\item{$C_{et, loc}$ the category of all  $G \in C_{et}$ with $G^D \in C_{loc}$,}
\item{$C_{et, et}$ the category of all  $G \in C_{loc}$ with $G^D \in C_{et}$.}
\end{itemize}

Hence the category $C$ has the following decomposition.
\[
C=C_{et, et}\times C_{et, loc} \times C_{loc, et} \times C_{loc, loc}
\]
Each category is abelian and hence $C$ is abelian itself.

\begin{proposition}
There are following equivalences:
\begin{description}
\item[a]{$G \in C_{et, et}$ if and only if ${\rm Frob_G}$ and ${\rm Ver}_G$ are isomorphisms, }
\item[b]{$G \in C_{et, loc}$ if and only if ${\rm Frob_G}$ is isomorphism and ${\rm Ver}_G$ is nilpotent,}
\item[c]{$G \in C_{loc, et}$ if and only if ${\rm Frob_G}$ is nilpotent  and ${\rm Ver}_G$ is isomorphism,}
\item[d]{$G \in C_{loc, loc}$ if and only if ${\rm Frob_G}$ and ${\rm Ver}_G$ are both nilpotents.}
\end{description}
\end{proposition} 

In this paper we focus on a group scheme $A[p]$, where $A$ is an Abelian variety defined over a finite field $\FF_{q}$, with $p={\rm char}(\FF_{q})$. The group scheme $A[p]$ belongs to the category $C$. Furthermore, due to the proposition above and  the relations 
${\rm Frob}\cdot{\rm Ver}={\rm Ver}\cdot{\rm Frob}=p$ and the fact that $p$ is nilpotent, it follows that
$$
A[p] \in C_{et, loc} \times C_{loc, et} \times C_{loc, loc}.
$$
\end{subsection}
\begin{subsection}{Dieudonn\'e functor}
Let $\mathbb{W}_{n}={\rm Spec}(k[X_1, \ldots, X_n])$ with an addition law given by the Witt polynomials.
Then $\mathbb{W}_{n}$ is called the $n-$th Witt group scheme. Denote
\begin{displaymath}
\mathbb{W}_{\infty}=\lim_{\leftarrow} \mathbb{W}_n
\end{displaymath}

Let $\mathbb{E}$ be the ring of "noncommutative polynomials" over $\mathbb{W}(k)$ in two variables $\mathcal{F}$ and $\mathcal{V}$ satisfying the following relations
\begin{itemize}
\item{$\mathcal{F} \, a=\sigma(a)\, \mathcal{F}$,}
\item{$\mathcal{V} \,\sigma(a)=a\, \mathcal{F}$,}
\item{$\mathcal{F}\mathcal{V}=\mathcal{V}\,\mathcal{F}=p$,}
\end{itemize}
where $\sigma$ is the Frobenius automorphism of $\mathbb{W}(k)$ and $a \in \mathbb{W}(k),$ i.\ e.\  $
\mathbb{E}=\mathbb{W}(k)[\mathcal{F}, \, \mathcal{V}]/ \sim
$

There exits the covariant Dieudonn\'e functor

\[
\mathbb{D}: \left\{
\begin{array}{ccl}
C_{et,loc}\times C_{loc, et} \times C_{loc, loc} &\rightarrow & 
\left\{
\begin{array}{c}
\mbox{left}\, 
\mathbb{E-}\mbox{modules of finite length}\\

 \end{array}\right\} \\
 &&\\
 G& \mapsto & {\rm Hom}(G, \mathbb{W}_{\infty})=\lim_{\to} {\rm Hom}(G, \mathbb{W}_{n}).
\end{array}
\right.
\]
 gives an equivalence of the category  $C_{et,loc}\times C_{loc, et} \times C_{loc, loc}$
and the category of the Dieudonn\'e modules.

\begin{remark}
In case of Abelian varieties,
$\mathbb{D}({\rm A}[p])$ is a vector space over $k$ of dimension  $2\,{\rm dim}(A)$ with 
actions $\mathcal{F}$ and $\mathcal{V}$.
\end{remark}

\end{subsection}
%------------------------------------------------------------------------------------
\begin{subsection}{ Kraft diagrams}
We will write ${\rm C}(1)_{k}$ for the category of finite commutative $k-$group schemes which are killed by $p$, where $k$ is a perfect algebraically closed field of characteristic  $p>0$.  The Dieudonn\'e functor shows that 
the full subcategory ${\rm C}(1)_{k}$ of $C$ is equivalent to the category of triples $({\rm M},\, {\mathcal F},\, {\mathcal V})$ where
\begin{itemize}
\item{${\rm M}$ is a finite dimensional $k$--vector space,}
\item{${\mathcal V}: {\rm M} \rightarrow {\rm M}$ is a ${\rm Frob}_{k}$-linear map,}
\item{${\mathcal F}: {\rm M} \rightarrow {\rm M}$ is a ${\rm Frob}_{k}^{-1}$-linear map,}
\end{itemize}
such that  ${\mathcal F}{\mathcal V}={\mathcal V}{\mathcal F}=0.$

\begin{definition}
A finite locally free commutative group scheme $G$ over a scheme $S$ is called 
a {\bf truncated Barsotti-Tate group of level $1$} or  a {\bf ${\rm BT}_{1}$} group scheme if it is killed by $p$ and ${\rm Ker}(F_{G})={\rm Im}(V_{G})$
\end{definition}

In terms of exact sequences, a group ${\rm G} \in {\rm C}(1)_{k}$ is a {\bf ${\rm BT}_{1}$} group scheme if and only if  the sequence 
$$
 \xymatrix{{\rm G}   \ar[r]^{\mathcal{F}_{\rm G}} &{\rm G}^{(p)} \ar[r]^{\mathcal{V}_{\rm G}} & {\rm G} }
$$
is exact. On the Dieudonn\'e module this means that
${\rm Ker}(\mathcal{V}) ={\rm Im}(\mathcal{F})$ and ${\rm Ker}(\mathcal{F}) ={\rm Im}(\mathcal{V}).$

In the unpublished manuscript \cite{Kraft}, Kraft  showes that 
there is a normal form of an object of ${\rm C}(1)_{k}$. The normal form is distinguished into two types of group schemes, linear and circular types and the group scheme $A[p]$ corresponds only to circular type.

\begin{subsubsection}{Circular type}
\begin{definition}
A {\bf circular word} is a finite ordered set of symbols $\mathcal{F}$ and $\mathcal{V}$:
$$
w=L_1 \ldots L_{t}, \qquad L_i \in \{\mathcal{F},\, \mathcal{V}\}.
$$
\end{definition}

We say that two words $w_1$ and $w_2$ are equivalent if there is a cyclic permutation which
transfers one word to another. In other words, the class of a word is given by 
$[L_1 \ldots L_{t}]=[L_2 \ldots L_{t}\, L_1]=\ldots=[L_t\, L_1 \ldots L_{t-1}].$

For a give word $w$ one can construct a finite group scheme $G_{w}$ over $k$ defined by
the $k-$vector space
\[
\mathbb{D}(G_{w})=\sum_{i=1}^{t} k\,z_{i},
\]
with structure of a Diuedonn\'e module given by

\begin{tabular}{l c}
if $L_i=\mathcal{F}$  & then $\mathcal{F} z_{i}:=z_{i+1}$ and $ \mathcal{V}\, z_{i+1}=0$,\\
if $L_i=\mathcal{V}$  & then $\mathcal{V} z_{i+1}:=z_{i}$ and $ \mathcal{F}\, z_{i}=0$.\\
\end{tabular}

There exits a visualization of it by  a circular graph $\Gamma_{w}$ where all edges are labeled by ${\mathcal F}$ and ${\mathcal V}$ pointing clockwise if $L_{i}=\mathcal{F}$
and counterclockwise otherwise.

 A word $w$ is called indecomposible if the corresponding graph $\Gamma_{w}$  is not periodic (i.\ e.\ the diagram is not invariant under non-trivial rotation). Two circular words equivalent if and only if the corresponding circular diagrams differ by a rotation. Since notions of  circular word and  circular diagram are equivalent we write  $G_{\Gamma}$ for a group scheme corresponding to a circular diagram $\Gamma$.
\end{subsubsection}
%\begin{subsubsection}{Cartier duality}
%The Cartier duality of $G_{\Gamma}$ is isomorphic to $G_{\hat{\Gamma}}$, where $\hat{\Gamma}$ is the diagram obtained from ${\Gamma}$ by changing all ${\mathcal F}-$arrow into ${\mathcal V}-$arrow and vice versa.

%\begin{proposition}\label{duality}
%The group scheme $G_{\Gamma}$ is a $BT_1$ if and only if $\Gamma$ is a circular diagram.
%\end{proposition}

%As a direct consequence of  Proposition~\ref{duality}  it follows that 
%a group scheme $G_{\Gamma}$ is a symmetric $BT_1$ if and only if $\Gamma$ is a circular diagram and $\Gamma=\hat{\Gamma}$. More

%\begin{definition}
%A  ${\rm BT}_1$ group scheme $G$ is called a symmetric  ${\rm BT}_1$ group scheme if $G=G^{D}$.
%\end{definition}
%\end{subsubsection}

The category $C(1)_{k}$ is abelian and all objects have finite length hence  ever object from  $C(1)_{k}$ is a direct sum of indecomposable objects. Up to isomorphism and permutation of the factors this decomposition is unique. The following theorem can found in  \cite{Kraft}.

\begin{theorem}
\begin{enumerate}
\item{A circular word $w$ defines a ${\rm BT}_1$ group scheme  $G_{w}$, and $w$ is indecomposible if and only if $G_{w}$ is indecomposible.}
\item{For any indecomposible ${\rm BT}_1$ group scheme $G$ there exits an indecomposible  word $w$ such that $G \cong G_{w}$. }
\end{enumerate}
 \end{theorem}
\end{subsection}

\begin{subsection}{The classification up to dimension $3$.}
Let $G$ be a finite commutative group scheme over  an algebraically closed field $\bar{\FF}_p$.
Then the classification of such groups schemes can be given either by the Kraft diagrams or the Dieudonne modules.
Here we produce a part of this classification for our needs.

%a symmetric truncated Barsotti-Tate group scheme over an algebraically closed field $\bar{\FF}_p$.
%Then there exits a classification of such group schemes via Kraft diagrams and Dieudonne modules of $BT_{1}$ group schemes.

\begin{tabular}{|c|c | c | c |}
\hline
Group scheme &Dieudonn\'e module  & Kraft diagram & ($a, \, p\mbox{-rank}$)\\
\hline
$\alpha_p$ & $ \EE/(F,V)_{l}$ & 

$$
    \xymatrix{\bullet }
$$
&($1$, $0$)\\
\hline
$\mu_p$ & $ \EE/(1-F,V)_{l}$ &

$$
    \xymatrix{\bullet  \ar@(ur,dr)[]^{\mathcal{F}}\\
    \\}
$$

&($0$, $1$)\\
\hline
$\ZZ/p\ZZ$ & $\EE/(F,1-V)_{l}$ &

$$
    \xymatrix{\bullet  \ar@(dl,ul)[]_{\mathcal{V}}\\
    \\}
$$

 &($0$, $0$)\\
\hline
${\rm I}_{1,1}$ & $\EE/(F+V)_{l}$ &
$$
    \xymatrix{
        \bullet \ar@/^/[r]^{\mathcal{F}}
         \ar@/_/@{>}[r]_{\mathcal{V}} &
        \bullet }
$$
& ($1$, $0$)\\
\hline
${\rm I}_{2,1}$ & $\EE/(F^2+V^2)_{l}$ &

$$
    \xymatrix{
        \bullet   & \bullet \ar@/_/[l]^{\mathcal{F}} \ar@/^/[d]^{\mathcal{F}} \\
        \bullet \ar@/_/[r]_{\mathcal{V}}  \ar@/^/[u]_{\mathcal{V}}     & \bullet}
$$

& ($1$, $0$)\\
\hline
${\rm I}_{3,1}$ & $\EE/(F^3+V^3)_{l}$ &
$$
    \xymatrix{
        &\bullet \ar@/^/[dr]^{\mathcal{F}} \ar@/_/[dl]^{\mathcal{V}} &\\
        \bullet\ar@/_/[d]^{\mathcal{V}}  && \bullet \ar@/^/[d]^{\mathcal{F}}\\
        \bullet\ar@/_/[dr]^{\mathcal{V}}  && \bullet \ar@/^/[dl]^{\mathcal{F}}\\
        &\bullet&        }
$$

& ($1$, $0$)\\
\hline
${\rm I}_{3,2}$ & $\EE/(F^2-V)_{l}\oplus \EE/(V^2-F)_{l}$ & 
$$
    \xymatrix{
    &\bullet \ar@/_/[dl]^{\mathcal{V}} \ar@/^/[d]^{\mathcal{F}} &   &   \bullet \ar@/_/[dl]^{\mathcal{F}} \ar@/^/[d]^{\mathcal{V}} \\
    \bullet & \ar@/^/[l]^{\mathcal{F}} \bullet&  \bullet & \ar@/^/[l]^{\mathcal{V}} \bullet }
$$

&($2$, $0$)\\
\hline
\end{tabular}
\end{subsection}
\end{section}

%------------------------------------------------------------------------------------
\begin{section}{groups schemes arising from ideals}

In this section we recall some basic properties of kernel ideals as group schemes, for details we refer to \cite{Waterhouse}, \cite{Shimura_Taniyma} (pp. 50--66) and \cite{Milne_Complex_Multiplication}(section~6). At the end of this section we give a proof of the Deuring reduction theorem in accordance with the theory of kernel ideals. 

Let $A$ be an abelian scheme over a scheme $S$ with complex multiplication by 
the full ring of integers $\cO_{K}$ of  a CM field $K$. For each isogeny $\alpha \in R$ the kernel $\Ker(\alpha)$ is a finite flat group scheme over $S$ denoted by $A[\alpha]$. For any 
ideal $I \subset \cO_{K}$ the intersection 
\[
A[I]= \cap_{\alpha \in I} A[\alpha]
\]
is a finite flat group scheme over $S$. 

\begin{proposition}\label{kernels}
Let $I$ and $J$ be ideals of $\cO_{K}$. Then the following hold
\begin{itemize}
\item{the order of group scheme $A[I]$ is ${\rm Nm}_{K/\QQ}(I),$}
\item{$A[I+J]=A[I] \cap A[J],$}
\item{$A[I\cap J]=A[I] + A[J],$}
\item{furthermore, $A[I \cdot J]=A[I] \oplus A[J]$, if $I$ and $J$ are coprime ideals,}
\item{$A[I]^{D}\cong A[I^{c}]$, where $A[I]^{D}$ is the Cartier dual.}
\end{itemize}
\end{proposition}

The next lemma gives a sufficient condition when  $A[p]$ is of local-local type.
\begin{lemma}\label{local-local}
Let $A$ be an Abelian variety over a finite field $\FF_{p^n}$ with complex multiplication by 
the full ring of integers $\cO_{K}$ of  a CM field $K$. Assume that 
\[
p\cO_{K}=\mathcal{P}_1 \ldots \mathcal{P}_{m}
\]
is a decomposition into prime ideals (not necessarily distinct) such that each prime ideal $\mathcal{P}$ is self-conjugated. Then the group scheme $A[p]$ is of local-local type.
\end{lemma}
\begin{proof}

Let $j:K \hookrightarrow {\rm End}(A) \otimes \QQ$ and $C$ be the center of ${\rm End}(A) \otimes \QQ$. Then $j^{-1}(C) \subset K $ and hence $j^{-1}(C \cap {\rm End}(A)) \subset \cO_{K}$. Therefore there is an element $\pi \in \cO_{K}$ corresponding to the relative Frobenius 
$F: A \to A^{(p^n)}=A$ (the $n-$the power of the Frobenius ${\rm Fr}: A \to A^{(p)}$), then
$\bar{\pi} \in \cO_{K}$ corresponds to the $n-$th power of the Verschiebung 
${\rm Ver}^{n}:A^{(p^n)}=A \to A$ and  $\pi \bar{\pi}=p^n$. Hence
$(\pi) (\bar{\pi})=(p)^n=\mathcal{P}_1^n \ldots \mathcal{P}_{m}^n$. Due to the uniqueness of the decomposition we conclude that $(\pi)=(\bar{\pi})$, therefore the Frobenius
and the Verschiebung  are nilpotents on $A[p]$ and hence
$A[p]$ is of local-local type.
\end{proof}

The following proposition is adopted version of Theorem 3.13 and  its proof from~\cite{Waterhouse}.
\begin{proposition}\label{maximal order}
Let $A$ be an Abelian variety over $\bar{\FF}_p$ with complex multiplication by 
a CM field $K$, i.\ e.\ $j:K \hookrightarrow {\rm End}(A) \otimes \QQ $. Then there exits an Abelian variety
$B$ over $\bar{\FF}_p$ with complex multiplication by the full ring of integer $\cO_{K}$ of $K$ which is isogenous to $A$.
\end{proposition}
\begin{proof}
Since  $\cD = j^{-1}({\rm End}(A) \cap j(K))$ is  a lattice of $K$ there is a positive integer $N$ such that 
$N \cO_{K} \subset \cD$. Let $I=\cD \cdot N\cO_{K}$ be an ideal of $\cD$ and $B$ be a quotient abelian variety $A$ by a group scheme $A[I]={\rm Ker}(I)$ and hence it is isogenous to $A$. Then endomorphism ring of $B$ contains $\cO_{K}$ and hence  $\cO_{K} \hookrightarrow {\rm End}(B)$.
\end{proof}

Now we can prove the Deuring theorem in terms of group schemes.
Namely, under the conditions of the Theorem \ref{Deuring}, we have that there is an elliptic curve ${E}^{\prime}$ with complex multiplication by the full ring of integer $\cO_{K}$ which is isogenous to $E:=(\cE \, \mbox{mod}\, \mathcal{P})$  over $\bar{\FF}_p$, by the Proposition~\ref{maximal order}, and hence  ${E}^{\prime}[p] \cong E[p]$ (note that it does not hold for higher dimension). 

If $p\cO_{K}=\mathcal{P}\mathcal{P}^c$ then 
$E[p] \cong{E}^{\prime}[p] \cong {E}^{\prime}[\mathcal{P}] \times {E}^{\prime}[\mathcal{P}]^{D}$, by the Proposition~\ref{kernels}. Thus $E[p] \cong \ZZ/p\ZZ \times \mu_{p}$ by the classification.

In case of $p$ is inert or ramified it follows that $E[p]$ is of local-local type, by 
the Lemma~\ref{local-local}. Thus $E[p] \cong {\rm I}_{1, 1}$.

In the next section we give more consistent approach to the problem of reduction of Abelian varieties.
\end{section}

%----------------------------------------------------------------------------------------------

\begin{section}{kraft diagrams and representation theory}\label{KraftRT}

In this section we show the decomposition group and a choice of a CM type  completely
determine  the ${\rm BT}_1$ group scheme of an Abelian variety that arises after the reduction of a CM Abelian variety at a place of good reduction. Here we develop an explicit representation theory 
on the base of the approach of \cite{Ekedahl}~(section~2).

Let $\mathcal{A}$ be an abelian  scheme of relative dimension $g$ over ${\rm Spec}(\cO_{L})$,
where $\cO_{L}$ is the full ring of integers of a number field $L$. Assume that $\mathcal{A}$
has a complex multiplication by the full ring of integers of a CM field $K$ and  $K \subset L$. Suppose that $\mathcal{A}$ has a good reduction 
at a prime ideal $\mathcal{B} \subset \cO_{L}$ and a number prime $p \in \mathcal{B} \cap \ZZ$ is unramified in $K$.

Let $p\cO_{K}=P_{1} \ldots P_{m}$ be a decomposition into distinct prime ideals.
Let $\tilde{K}$ be a Galois closure of $K$ over $\QQ$  and $p\cO_{K}= \tilde{P}_{1} \ldots \tilde{P}_{l}$ be a decomposition into  prime ideals in $\tilde{K}$ (note, that $p$ is also unramified in the Galois closure $\tilde{K}$, since it is a composite of all embeddings of $K$ into fixed algebraic closure of $\QQ$).
%Let $\mathcal{P}$ be a prime of $\cO_{L}$ and $p$ be a prime number in $\mathcal{P} \cap \ZZ$.

%%%%%%%%%%%%
The ring $\cO_{K}$ is a free $\ZZ-$module and $\cO_{K} \otimes_\ZZ \QQ \cong K.$ The semisimple $\bar{\QQ}$-algebra $\cO_{K} \otimes_\ZZ \bar{\QQ}$ can be decomposed into  irreducible components
\begin{displaymath}
\cO_{K} \otimes_\ZZ \bar{\QQ}\cong\prod_{\alpha \in {\rm Hom}(K, \bar{\QQ})} \bar{\QQ}
\end{displaymath}
and hence it induces a decomposition into simple $\cO_{K} \otimes_\ZZ \bar{\FF}_p$-modules
\begin{displaymath}
\cO_{K} \otimes_\ZZ \bar{\FF}_p\cong\prod_{\alpha \in {\rm Hom}(K, \bar{\QQ})} \bar{\FF}_p,
\end{displaymath}
since $p$ is unramified in $K$. So the representation on ${\rm H^1_{DR}}(A/ \bar{\FF}_p)$ of the semi-smiple algebra
$\cO_{K}\otimes_\ZZ \bar{\FF}_p$   is a direct sum of irreducible representations
\begin{displaymath}
{\rm H^1_{DR}}(A/ \bar{\FF}_p) \cong \oplus_{\alpha \in {\rm Hom}(K, \CC)} V_{\alpha},
\end{displaymath}
where $V_{\alpha}$ is an irreducible $\cO_{K} \otimes_\ZZ \bar{\FF}_p$-module.  

%%%%%%%%%
In order to get an explicit description of the Kraft diagram we reformulate this decomposition in terms of Galois actions. 
Let $G$ be a Galois group ${\rm Gal}(\tilde{K}/ \QQ)$, $\Delta:={\rm Gal}(\tilde{K}/K)$, $\tilde{P}=\mathcal{B} \cap \cO_{\tilde{K}}$  and $\sigma$ be a generator of the decomposition group of a prime ideal $\tilde{P}$ in $\tilde{K}$ . 
Then the decomposition of $\cO_{K} \otimes_\ZZ \bar{\FF}_p$ can be written as following
\begin{displaymath}
\cO_{K}/p\cO_{K} \otimes_\ZZ \bar{\FF}_p\cong
\oplus_{i=1}^{m} F_{P_i} \otimes_{\FF_p} \bar{\FF}_p \cong
\oplus_{i=1}^{m}\left( \oplus_{\alpha \in {\rm Hom}(F_{P_i}, \bar{\FF}_p)}  \bar{\FF}_p\right)
 \cong\prod_{\alpha \in {\rm Hom}(K, \bar{\QQ})} \bar{\FF}_p,
\end{displaymath}
where $F_{P_i}$ is the residue field of a prime ideal $P_i$. The last isomorphism comes 
from the fact that the embeddings $F_{P_i} \to \bar{\FF}_p$ one-to-one correspond to embeddings $K \to \bar{\QQ}$ (since $p$ is unramified).

Let's fix a prime ideal  $\tilde{P}:=\tilde{P}_i$  and an isomorphism $\cO_{\tilde{K}}/ \tilde{P} \cong \FF_{q} \subset \bar{\FF}_{p}$. Then each $\alpha \in G$ induces  an 
embedding $\FF_{q}$ into  $\bar{\FF}_{p}$, by sending $(a \,{\rm mod } \, \tilde{P}) \mapsto (\alpha(a) \,{\rm mod } \, \tilde{P})$. 
Then we have the following decomposition of ${\rm H^1_{DR}}(A/ \bar{\FF}_p)$
into $2g$ one dimensional eigenspaces

$$
{\rm H^1_{DR}}(A/ \bar{\FF}_p)=\bigoplus_{\alpha \in G \setminus \Delta} {\rm V}_{\alpha},
$$

where $\alpha$ runs all conjugate classes of $G$ by action of $\Delta$ on the right and
$$
{\rm V}_{\alpha}=\{v \in {\rm H^1_{DR}}(A/ \bar{\FF}_p) \, | \, a(v)= (\alpha(a) \,{\rm mod } \, \tilde{P}) v\qquad \mbox{for any}\, a \in \cO_{K}    \}.
$$
Since all $V_{\alpha}$ are isomorphic to each other it follows that ${\rm dim}_{\bar{\FF}_p}V_{\alpha}=1$ for each $\alpha \in G \setminus \Delta$.

%%%%%%%%%%%

The action of $\cO_K$ on ${\rm H^1_{DR}}(A/ \bar{\FF}_p)$ is imposed by a fixed isomorphism $\cO_{K} \cong {\rm End}(A)$. The Frobenius ${\rm Fr}$ on ${\rm H^1_{DR}}(A/ \bar{\FF}_p)$ is $p$-linear and it commutes with the  $\cO_K$-action and, hence 
the ${\rm Fr}(V_{\alpha})$ is $\cO_{K} \otimes_\ZZ \bar{\FF}_p$-module and it can be decomposed into irreducible components. Each irreducible component is one dimensional $\bar{\FF}_p$-vectors space and the dimension  $\bar{\FF}_p$-vector space ${\rm Fr}(V_{\alpha})$ is less or equal to $1$. Therefore for each
irreducible component $V_{\alpha}$ of ${\rm H^1_{DR}}(A/ \bar{\FF}_p)$ there exits an irreducible component $V_{\beta}$ such that ${\rm Fr}(V_{\alpha}) \subset V_{\beta}$. Moreover, $\beta$ corresponds to the class of $\sigma \alpha$ in $G \setminus \Delta$
(since $\sigma(a) \equiv a^p\, \mbox{mod}\, \tilde{P}$). In other words,
$$
{\rm Fr}:{\rm V}_{\alpha} \rightarrow {\rm V}_{\sigma \alpha}.
$$
%%%%%%%%%%

Let's denote the fibre product $(\mathcal{A} \, {\rm mod}\, \mathcal{B}) \times \bar{\FF}_{p}$ by $A$. The set $S$ of  the isomorphism classes of irreducible factors of ${\rm H^1_{DR}}(A/ \bar{\FF}_p)$ as $\cO_K\otimes_{\ZZ} \bar{\FF}_p$-module can be identify with
a set $\Hom(K, \bar{\QQ})$ via identification of the isomorphism classes of  irreducible representations of $\cO_K\otimes_{\ZZ} \bar{\FF}_p$ in $\bar{\FF}_p$ and the set $S$.
Due the canonical exact sequence 
$$
0 \rightarrow {\rm H^0}(A, \Omega^1_{A})\rightarrow {\rm H^1_{DR}}(A/ \bar{\FF}_p) \rightarrow {\rm H^1}(A, \cO_A)\rightarrow 0
$$
$S$ is disjoint union of $S^0$ and $S^1$, where $S^0$ and $S^1$ are the classes of irreducible representations occurring in 
${\rm H^0}(A, \Omega^1_{A})$ and ${\rm H^1}(A, \cO_A)$, respectively. 
\par
In order to describe the action of the Frobenius and the Verschiebung morphisms on the Dieudonn\'e module ${\mathbb D}(A[p])$ we apply the equivalence of the categories of the first de Rham cohomology and the Dieudonn\'e modules, see \cite{Oda} (with exception that we use here the covariant Dieudonn\'e theory). There is an exact sequence of group schemes
$$
0 \rightarrow 
A[{\rm Ver}]
\rightarrow 
A[p]
\xrightarrow{\rm Ver}
A[{\rm Fr}]
\rightarrow 0,
$$
which yields an exact sequence of the Dieudonn\'e modules
$$
0 \rightarrow 
{\mathbb D}(A[{\rm Ver}])
\rightarrow 
{\mathbb D}(A[p])
\xrightarrow{\mathcal F} 
{\mathbb D}(A[{\rm Fr}])
\rightarrow 0.
$$
Applying the equivalence of the categories we obtain the following commutative diagram
(for example, see \cite{Goren_book}, pp.\ 245--248)
\begin{displaymath}
\xymatrix{
        0  \ar[r] & {\rm H^0}(A, \Omega^1_{A}) \ar[d]^{\cong} \ar[r] &
        {\rm H^1_{DR}}(A/ \bar{\FF}_p) \ar[d]^{\cong} \ar[r]&  {\rm H^1}(A, \cO_A)\ar[d]^{\cong} \ar[r]& 0 \\
       0  \ar[r] & {\mathbb D}(A[{\rm Ver}]) \ar[r]&  {\mathbb D}(A[p]) \ar[r]& {\mathbb D}(A[{\rm Fr}]) \otimes_{\rm Fr} \bar{\FF}_p\ar[r]& 0.}
\end{displaymath}

%%%%%%%%%%%%%%%%%%%%%%
Let $S^{1}$ be a CM type of $\mathcal{A}$ and $S^{0}$ be the conjugation of the CM type. 
On the base of these data we can construct the Kraft diagram
for ${\rm BT}_1$ of an Abelian variety $A$ over  $\bar{\FF}_{p}$ 
(and, as a consequence, we obtain $p$-rank and $a$-number of $A$).
%%%%%%%%%%%%%%%

Now we utilize the fact that ${\rm H^0}(A, \Omega^1_{A}) \cong {\mathbb D}(A[{\rm Ver}])$, it implies on the one hand that ${\rm Frob}({\rm V}_{\alpha})=0$  and on the other hand an ismorphism ${\rm Ver}:{\rm V}_{\sigma \alpha} \tilde{\rightarrow} {\rm V}_{\alpha}$ for each $\alpha \in S^{0}$ (here we identify $S^0$ with corresponding subset of classes in $G \setminus \Delta$, under the identification $S$ with $G \setminus \Delta$).

Therefore on the base of CM type $S^1$, group $\Delta$, a  Galois group $G={\rm Gal}(\tilde{K}/\QQ)$ and a decomposition group of a prime ideal in $\tilde{K}$ over $p$, one can construct a Kraft diagram $\Gamma$. Vertices of $\Gamma$ are classes $G / \Delta$ and there is an arrow between a  class $[\alpha] \in G / \Delta$
and $[\sigma \alpha]$, where $\sigma$ is a generator of  the decomposition group. The arrow  $[\alpha]\to[\sigma \alpha]$ is labeled by $\mathcal{V}$  if $[\alpha] \in S^1$ and  $[\sigma \alpha]\to[\alpha]$ is labeled by $\mathcal{F}$ otherwise.
\end{section}
%%%%%%%%%%%%%%%%%%%%%%%%%%%%%%%%%%%%%%%
\begin{section}{abelian schemes of relative dimension $1, \, 2$ and $3$}

In this section we prove our main theorem~\ref{Main_Theorem} stated in Section~\ref{introduction}. 

For given a Galois group $G:={\rm Gal}(\tilde{K}/\QQ)$ of  a Galois closure $\tilde{K}$ of $K/\QQ$, $\Delta={\rm Gal}(\tilde{K}/K)$, the complex involution $\iota$ and a generator $\sigma$ of the decomposition ${\rm D}(\mathcal{P}) \subset G$, one can find a decomposition of prime $p$ into prime ideals in $K$ which corresponds to the group $\Delta \setminus {\rm Gal}(\tilde{K}/\QQ) /< \sigma>$. If additionally, a CM type is given then it is possible to construct a corresponding Kraft diagram as it is explained in Section~\ref{KraftRT}.

\begin{subsection}{Elliptic Schemes}
In this subsection we reproduce the Deuring result on the base of previous section.

Let $\mathcal{E}$ be an elliptic scheme over ${\rm Spec}(\cO_{L})$, where $\cO_{L}$ is the ring
of integers of a number field $L$. Assume that $\mathcal{E}$ has a complex multiplication 
by  the full ring of integers $\cO_{K}$ of a quadratic imaginary field $K$ and $K \subset L$. For given prime $\mathcal{B}$ in $\cO_{L}$, the reduction $E =\mathcal{E} \, {\rm mod} \, \mathcal{B}$ at a place of good reduction $\mathcal{B}$  is defined over $\cO_{L}/\mathcal{B}$.

The Galois group ${\rm Gal}(K/\QQ)=\{1, x\} \cong \ZZ/2\ZZ$ and, hence we may assume that CM-type $S^{1}=\{1\}$ and  $S^{0}=\{x \}$.

If  $p$ splits completely, i.\ e.\  $p\cO_{K}=\mathcal{P}_{1}\mathcal{P}_{2}$
then the decomposition group  is trivial , ${\rm D}(\mathcal{P}_{1})=\{1\}$, and hence a generator of the decomposition group  $\sigma=1$.  The action of $\sigma$  is the following:
\[
\sigma \cdot1=1 \qquad \sigma \cdot x=x.
\] 
In terms of circular words and diagrams, we have  two circular words $\omega_{1}=[\mathcal{F}]$ and $\omega_{2}=[\mathcal{V}]$ or the Kaft diagrams

\begin{displaymath}
    \xymatrix{\bullet  \ar@(ur,dr)[]_{\mathcal{F}}\\
    \\}\qquad
    \xymatrix{\bullet  \ar@(dl,ul)[]_{\mathcal{V}}\\
    \\}
\end{displaymath}

Therefore $E[p] \times \bar{\FF}_{p} \cong \ZZ/ p\ZZ \times \mu_{p}$.

If $p$ is inert, the decomposition group $D(\mathcal{P})=<x>$ and the action of a generator is the follwoing
$$
\sigma  \cdot 1=x \qquad \sigma \cdot x=1,
$$ 
hence  a circular word is  $\omega=[\mathcal{F},\mathcal{V}]$,  which corresponds to the Kraft diagram 

$$
    
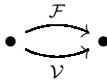
\xymatrix{
        \bullet \ar@/^/[r]^{\mathcal{F}}
         \ar@/_/@{>}[r]_{\mathcal{V}} &
        \bullet }
$$

So $E[p] \times \bar{\FF}_{p} \cong I_{1,1}$.

Combining  all result above we have the following table.

\begin{tabular}{|c|c | c | c | c |}
\hline
ideal decomposition &$BT_{1}$ group scheme & $p$-rank & $a$-number \\
\hline
$PP^{c}$ & $ \ZZ/p \times \mu_{p}$ & $1$ & $0$ \\
\hline
$P$ & $I_{1,1}$& $0$&$1$ \\
\hline
\end{tabular}
\end{subsection}

\begin{subsection}{Abelian Surfaces}

Let $\mathcal{A}$ be a abelian scheme of relative dimension $2$ over ${\rm Spec}(\cO_{L})$, where $\cO_{L}$ is the ring
of integers of a number field $L$. Assume that $\mathcal{A}$ has a complex multiplication 
by $\cO_{K}$, the full ring of integers of a CM field $K$, and $K \subset L$. For given prime $\mathcal{B}$ in $\cO_{L}$ of good reduction, the reduction $A =\mathcal{A} \, {\rm mod} \, \mathcal{B}$ is defined over $\cO_{L}/\mathcal{B}$.

There are two cases, namely, $K/\QQ$ is Galois and it is not Galois.

\begin{subsubsection}{  $K/\QQ$ is Galois}
The Galois case is simple, therefore we summarize all results together.

Let $K/\QQ$ be a Galois extension. Then a Galois group of $K/\QQ$ is cyclic, say ${\rm Gal}(K/\QQ)=\{1,\, x, \, x^{2}, \, x^{3} \}$, $\Delta=\{1\}$ and the complex conjugation corresponds to the element $x^{2}$. There is only one choice of  CM-type, $S^{1}=\{1, x \}$ and $S^{0}=\{x^{2}, x^{3} \}$.
Thus 

\begin{tabular}{|c|c|c|c|}
\hline
$\sigma$ & $p\cO_{K}$ & word/words & $A[p] \times \bar{\FF}_{p}$\\
\hline
$1$ & $\mathcal{P}_{1} \mathcal{P}_{2} \mathcal{P}^{c}_{1} \mathcal{P}^{c}_{2}$ & 
$[\mathcal{F}],\, [\mathcal{F}], [\mathcal{V}],\, [\mathcal{V}] $ & $(\ZZ/p\ZZ)^2\times \mu_{p}^2$\\
\hline
$x$ & $\mathcal{P}$ & $ [ \mathcal{F}\,  \mathcal{F}\,  \mathcal{V}\,  \mathcal{V}]$& ${\rm I}_{1,2}$\\
\hline
$x^2$ & $\mathcal{P}_1 \mathcal{P}_2$ & $[\mathcal{F}, \mathcal{V}], \, [\mathcal{F}, \mathcal{V}]$ & ${\rm I}_{1,1}^2$\\
\hline
\end{tabular}
\end{subsubsection}

\begin{subsubsection}{ $K/\QQ$ is non-Galois extension.}
If $K/\QQ$ is not a Galois extension then 
$G:={\rm Gal}(\tilde{K}/\QQ) \cong D_{4}=<x,y\, |\, y^{4}=x^{2}=xyxy=1>$, the
dihedral group of order 8. We always may assume that a primitive CM-type $S^{1}$ is $\{1,\, y\}$ and $S^{0}=\{y^{2},\, y^{3}\}$ (since $y^{2}$ is unique involution in the center of $D_{4}$).
Moreover, $\tilde{K}^{<x>}=K$  (i.\ e.\  $\Delta=\{1,\,x\}$) and hence
$S^{1}$ is identified with a subset of classes $\{\{1,\, x\}, \, \{y,\, yx\} \}$ 
(and  $S^{0}$ with $\{\{y^{2},\, y^{2}x\}, \, \{y^{3},\, y^{3}x\} \}$) in $G \setminus \Delta$.

Now we produce detail calculation in case of 
$\sigma=x$. The factor group  $$\Delta \setminus {\rm Gal}(\tilde{K}/\QQ) /< \sigma>$$ can be written as following
$$
[\{1,\, x\}],\, [\{y^2, \, xy^2\}], \, [\{y, \, yx, \, xy, \, y^3\}],
$$
and, hence
$$
p\cO_{K}=\mathcal{P}_{[\{1,\, x\}]} \mathcal{P}_{[\{y^2, \, xy^2\}]} 
\mathcal{P}_{[\{y, \, yx, \, xy, \, y^3\}]}=\mathcal{P}_1 \mathcal{P}_1^{c} \mathcal{P}_2,
$$
where $\mathcal{P}_2$ is self-conjugated.

On the other hand 
$$
\begin{array}{lc}
\sigma \{1,\, x\} \mapsto \{x,\, 1 \} &   \xymatrix{\bullet  \ar@(ur,dr)[]_{\mathcal{F}}  &  }\\
\sigma \{y,\,yx\} \mapsto \{xy,\, y^3 \} &   \xymatrix{\bullet  \ar[r]^{\mathcal{F}} & \bullet }\\
\sigma \{y^3,\, y^3x \} \mapsto \{xy^3,\, y \}&   \xymatrix{\bullet  & \bullet \ar[l]^{\mathcal{V}} }\\
\sigma \{y^2,\, y^2x \} \mapsto \{x y^2,\, y^2 \} &   \xymatrix{\bullet  \ar@(ur,dr)[]_{\mathcal{V}}  &  }\\
\end{array}
$$
so the Kraft diagram is
$$
\xymatrix{\bullet  \ar@(ur,dr)[]_{\mathcal{F}} &  }\,
\xymatrix{\bullet  \ar@(ur,dr)[]_{\mathcal{V}}  &  }\,
 \xymatrix{
        \bullet \ar@/^/[r]^{\mathcal{F}}
         \ar@/_/@{>}[r]_{\mathcal{V}} &
        \bullet }
$$
in other words $A[p] \times \bar{\FF}_{p} \cong \ZZ/p\ZZ \times \mu_{p} \times {\rm I}_{1,1}$.

%%%%%%%%%%%%%%%
Now we summarize all computations above in the following table

\begin{tabular}{|c|c|c|c|}
\hline
$\sigma$ & $p\cO_{K}$ & word/words & $A[p] \times \bar{\FF}_p$\\
\hline
$1$ & $\mathcal{P}_1\mathcal{P}_1^{c}\mathcal{P}_2\mathcal{P}_2^{c}$ &
$[\mathcal{F}], \,[\mathcal{F}], \, [\mathcal{V}], \, [\mathcal{V}]$ & $(\ZZ/p\ZZ \times \mu_{p})^2$\\
\hline
$xy$ & $\mathcal{P}\mathcal{P}^{c}$ &
$[\mathcal{F}], \,[\mathcal{F}], \, [\mathcal{V}], \, [\mathcal{V}]$ & $(\ZZ/p\ZZ \times \mu_{p})^2$\\
\hline
 $xy^3$ & $\mathcal{P}\mathcal{P}^{c}$ & $[\mathcal{F}\,\mathcal{V}],\, [\mathcal{F}\,\mathcal{V}]$ & ${\rm I}_{1,1}^2$\\
\hline
$x$ & $\mathcal{P}_1 \mathcal{P}_1^{c} \mathcal{P}_2$ & 
$[\mathcal{F}], \,[\mathcal{V}], \, [\mathcal{F} \,\mathcal{V}]$ 
& $\ZZ/p\ZZ \times \mu_{p} \times {\rm I}_{1,1}$\\
\hline
$y$ or $y^3$ & $\mathcal{P}$ & $[\mathcal{F} \,\mathcal{F} \, \mathcal{V} \,\mathcal{V}]$ &
${\rm I}_{2,1}$\\
\hline

\end{tabular}

%%%%%%%%%%%%
\end{subsubsection}

\begin{remark}
The same result was also obtained  by E.\ Goren in the paper~\cite{Goren}. In this paper the last case was solved on the base of a criterion on an Abelian surface  is superspecial.
\end{remark}
\end{subsection}

\begin{subsection}{Abelian Threefolds}
%%%%%%%%%%%%%%
Let $\mathcal{A}$ be a abelian scheme of relative dimension $3$ over ${\rm Spec}(\cO_{L})$, where $\cO_{L}$ is the ring
of integers for some number field $L$. Assume that $\mathcal{A}$ has a complex multiplication 
by $\cO_{K}$ the full ring of integers of a CM field $K$ and $K \subset L$. For  given prime $\mathcal{B}$ in $\cO_{L}$ of good reduction, the reduction $A =\mathcal{A} \, {\rm mod} \, \mathcal{B}$ is defined over $\cO_{L}/\mathcal{B}$.

According to the paper \cite{CM-Fields}, page~20
 $$
 {\rm Gal}(\tilde{K}/\QQ) \in  
 \{ \ZZ_{6}, \ZZ_{2}\times S_{3}, \ZZ_{2}^{3} 
 \rtimes \ZZ_{3},  (\ZZ_{2})^{3} \rtimes S_{3} \}.
 $$

\begin{subsubsection}{The Galois group is isomorphic to $(\ZZ/2\ZZ)^{3} \rtimes S_{3}$}

 Let's consider the case   $G:={\rm Gal}(\tilde{K}/\QQ) \cong (\ZZ/2\ZZ)^{3} \rtimes S_{3}$
 in details.  Other three case are less complicated and they can be work out in the similar manner.
 
 We start with some well known facts about the group 
 $$S_{3}=\langle s,\, t | s^{3}=t^{2}=1, \, ts=s^{2}t\rangle.$$ The center of $S_{3}$ is trivial, $\langle s \rangle$ is a unique normal subgroup and $\#<s>=3$, there are the following isomorphisms 
\begin{itemize}
\item{$S_{3} \rightarrow S_{3}$ given by $t \mapsto st$,}
\item{$S_{3} \rightarrow S_{3}$ given by $t \mapsto ts$,}
\end{itemize}
(roughly speaking, $t, \, st$ and $ts$ are interchangeable). 
 
A semi-direct product  $(\ZZ/2\ZZ)^{3} \rtimes S_{3}$ can be constructed via a homomorphism 
$S_{3} \to {\rm Aut}((\ZZ/2\ZZ)^{3})$ given by
$s \cdot (a_1,\, a_2,\, a_3)=(a_2, \, a_3,\, a_1)$ and
$t \cdot (a_1,\, a_2,\, a_3)=(a_2, \, a_1,\, a_3).$ Thus the group
 $(\ZZ/2\ZZ)^{3} \rtimes S_{3}$ is a set 
 $$
 \{(a_1,\, a_2,\, a_3;\, \alpha)| \qquad a_1, a_2, a_3 \in \ZZ/2\ZZ, \, \alpha \in S_3\}
 $$
with a group operation given by the following rule
$$
(a_1,\, a_2,\, a_3;\, \alpha) \cdot (b_1,\, b_2,\, b_3;\, \beta)=
((a_1,\, a_2,\, a_3)+\alpha \cdot (b_1,\, b_2,\, b_3) ;\, \alpha \cdot \beta).
$$ 

From the paper~\cite{CM-Fields}, it follows that the complex involution $\iota$
 corresponds to the element $(1,1,1;1)$ and ${\rm Gal}(\tilde{K}_{0}/\QQ) \cong S_3$.

On the base of all above, we always may assume that a subgroup $\Delta$ of  ${\rm Gal}(\tilde{K}/\QQ)$, fixing $K$, is generated by $(1,0,0;1),\, (1,0,0;1)$ and $(0,0,0;t)$. Then the factor group 
$$
G/ \Delta=
\left\{
\begin{array}{l}
\phi_1=[(0,0,0;1)],\,
\phi_2=[(0,0,1;1)]=\iota \phi_1,\,
\phi_3=[(0,0,0;st)],\, \\
\phi_4=[(0,0,0;ts)],\, 
\phi_5=[(1,1,1;st)]=\iota \phi_3,\,
\phi_6=[(1,1,1;ts)]=\iota \phi_4
\end{array}
\right\}
$$
Hence  the CM-types can be written in terms of elements of the factor group. There are four isomorphism classes of CM-types, namely:
\begin{description}
\item[(A)]{ $S_{1}:=\{   \phi_{1}, \, \phi_{2}, \, \phi_{4} \},$  }
\item[(B)]{$S_{1}:=\{   \iota\phi_{1}, \, \phi_{3}, \, \phi_{4} \},$ }
\item[(C)]{$S_{1}:=\{   \phi_{1}, \, \iota\phi_{3}, \, \phi_{4} \},$}
\item[(D)]{$S_{1}:=\{   \phi_{1}, \, \phi_{3}, \, \iota\phi_{4}\}.$}
\end{description}

Now we can carry out the similar computations as in the case of dimension two for each  CM-type,  and produce the similar table, replacing   the column "word/words" by  a column "CM-types".

\begin{table}[h]
\caption {Elements of order $2$} \label{order two} 
%%%%%%----------table---------%%%%%%%%%
\begin{tabular}{|c|c|c|c|}
\hline
$\sigma$ & $p\cO_{K}$ & CM-types & $A[p] \times \bar{\FF}_p$\\
\hline
%%%%%%%%%%
\multirow{2}{*}{$(0,0,1;ts)$} & 

\multirow{2}{*}{$\mathcal{P}_1\mathcal{P}_1^{c}\mathcal{P}_2\mathcal{P}_2^{c}$} &

(D), (B) &

$(\ZZ/p\ZZ \times \mu_{p})^3$\\ \cline{3-4}

&
&
(A), (C) &
$(\ZZ/p\ZZ \times \mu_{p})\times {\rm I}_{1, 1}^2$\\
\hline
%%%%%%%%%
\multirow{2}{*}{$(0,0,0;ts)$} & 

\multirow{2}{*}{$\mathcal{P}_1\mathcal{P}_1^{c}\mathcal{P}_2\mathcal{P}_2^{c}$} &

(A), (C) &

$(\ZZ/p\ZZ \times \mu_{p})^3$\\ \cline{3-4}

&
&
(B), (D) &
$(\ZZ/p\ZZ \times \mu_{p})\times {\rm I}_{1, 1}^2$\\
\hline
%%%%%%%%%
\multirow{2}{*}{$(0,0,1;t)$} & 

\multirow{2}{*}{$\mathcal{P}_1\mathcal{P}_1^{c}\mathcal{P}_2$} &

(A), (B) &

${\rm I}_{1,1}^3$\\ \cline{3-4}

&
&
(C), (D) &
$(\ZZ/p\ZZ \times \mu_{p})^2\times {\rm I}_{1, 1}$\\
\hline
%%%%%%%%%
\multirow{2}{*}{$(0,0,0;ts^2)$} & 

\multirow{2}{*}{$\mathcal{P}_1\mathcal{P}_1^{c}\mathcal{P}_2$} &

(A), (D) &

$(\ZZ/p\ZZ \times \mu_{p})^2\times {\rm I}_{1, 1}$\\ \cline{3-4}

&
&
(b), (C) &
${\rm I}_{1,1}^3$\\
\hline
%%%%%%%%%
\multirow{2}{*}{$(1,0,0;ts^2)$} & 

\multirow{2}{*}{$\mathcal{P}_1\mathcal{P}_1^{c}\mathcal{P}_2\mathcal{P}_2^c$} &

(A), (D) &

$(\ZZ/p\ZZ \times \mu_{p})^3$\\ \cline{3-4}

&
&
(B), (C) &
$(\ZZ/p\ZZ \times \mu_{p}) \times {\rm I}_{1,1}^2$\\
\hline
%%%%%%%%%
\multirow{2}{*}{$(0,1,0;ts)$} & 

\multirow{2}{*}{$\mathcal{P}_1\mathcal{P}_1^{c}\mathcal{P}_2$} &

(A), (C) &

$(\ZZ/p\ZZ \times \mu_{p})^2\times {\rm I}_{1,1}$\\ \cline{3-4}

&
&
(B), (D) &
$ {\rm I}_{1,1}^3$\\
\hline
%%%%%%%%%
\multirow{2}{*}{$(1,1,1;ts)$} & 

\multirow{2}{*}{$\mathcal{P}_1\mathcal{P}_1^{c}\mathcal{P}_2$} &

(B), (D) &

$(\ZZ/p\ZZ \times \mu_{p})^2\times {\rm I}_{1,1}$\\ \cline{3-4}

&
&
(A), (C) &
$ {\rm I}_{1,1}^3$\\
\hline
%%%%%%%%%
\multirow{2}{*}{$(1,1,0;t)$} & 

\multirow{2}{*}{$\mathcal{P}_1\mathcal{P}_1^{c}\mathcal{P}_2\mathcal{P}_2^c$} &

(A), (C) &

$(\ZZ/p\ZZ \times \mu_{p})^3$\\ \cline{3-4}

&
&
(B), (D) &
$(\ZZ/p\ZZ \times \mu_{p}) \times {\rm I}_{1,1}^2$\\
\hline
%%%%%%%%%
\multirow{2}{*}{$(1,1,1;ts^2)$} & 

\multirow{2}{*}{$\mathcal{P}_1\mathcal{P}_1^{c}\mathcal{P}_2\mathcal{P}_2^c$} &

(B), (C) &

$(\ZZ/p\ZZ \times \mu_{p})^3$\\ \cline{3-4}

&
&
(A), (D) &
$(\ZZ/p\ZZ \times \mu_{p}) \times {\rm I}_{1,1}^2$\\
\hline
%%%%%%%%%
\multirow{2}{*}{$(0,1,1;ts^2)$} & 

\multirow{2}{*}{$\mathcal{P}_1\mathcal{P}_1^{c}\mathcal{P}_2$} &

(B), (C) &

$(\ZZ/p\ZZ \times \mu_{p})^2\times {\rm I}_{1,1}$\\ \cline{3-4}

&
&
(A), (D) &
$ {\rm I}_{1,1}^3$\\
\hline
%%%%%%%%%
\multirow{2}{*}{$(0,0,0;t)$} & 

\multirow{2}{*}{$\mathcal{P}_1\mathcal{P}_1^{c}\mathcal{P}_2\mathcal{P}_2^c$} &

(D), (C) &

$(\ZZ/p\ZZ \times \mu_{p})^3$\\ \cline{3-4}

&
&
(A), (B) &
$(\ZZ/p\ZZ \times \mu_{p}) \times {\rm I}_{1,1}^2$\\
\hline
%%%%%%%%%
\multirow{2}{*}{$(0,0,0;ts^2)$} & 

\multirow{2}{*}{$\mathcal{P}_1\mathcal{P}_1^{c}\mathcal{P}_2\mathcal{P}_2^c$} &

(A), (D) &

$(\ZZ/p\ZZ \times \mu_{p})^3$\\ \cline{3-4}

&
&
(B), (C) &
$(\ZZ/p\ZZ \times \mu_{p}) \times {\rm I}_{1,1}^2$\\
\hline
%%%%%%%%%
\multirow{2}{*}{$(0,0,0;ts)$} & 

\multirow{2}{*}{$\mathcal{P}_1\mathcal{P}_1^{c}\mathcal{P}_2\mathcal{P}_2^c$} &

(A), (C) &

$(\ZZ/p\ZZ \times \mu_{p})^3$\\ \cline{3-4}

&
&
(B), (D) &
$(\ZZ/p\ZZ \times \mu_{p}) \times {\rm I}_{1,1}^2$\\
\hline
%%%%%%%%%
\multirow{2}{*}{$(1,1,1;t)$} & 

\multirow{2}{*}{$\mathcal{P}_1\mathcal{P}_1^{c}\mathcal{P}_2$} &

(A), (B) &

$(\ZZ/p\ZZ \times \mu_{p})^2\times {\rm I}_{1,1}$\\ \cline{3-4}

&
&
(C), (D) &
$ {\rm I}_{1,1}^3$\\
\hline
%%%%%%%%%
\multirow{2}{*}{$(1,1,1;st)$} & 

\multirow{2}{*}{$\mathcal{P}_1\mathcal{P}_1^{c}\mathcal{P}_2\mathcal{P}_2^c$} &

(B), (C) &

$(\ZZ/p\ZZ \times \mu_{p})^3$\\ \cline{3-4}

&
&
(A), (D) &
$(\ZZ/p\ZZ \times \mu_{p}) \times {\rm I}_{1,1}^2$\\
\hline
%%%%%%%%%
\multirow{2}{*}{$(1,1,1;ts)$} & 

\multirow{2}{*}{$\mathcal{P}_1\mathcal{P}_1^{c}\mathcal{P}_2$} &

(C), (D) &

$(\ZZ/p\ZZ \times \mu_{p})^2\times {\rm I}_{1,1}$\\ \cline{3-4}

&
&
(A), (B) &
$ {\rm I}_{1,1}^3$\\
\hline
%%%%%%%%%%%%%%%%%%%%%%%%%%%%%%%%
%%%%%%%%%%%%%%%%%%%%%%%%%%%%%%%%
$(1,1,1;1)$ & 

$\mathcal{P}_1\mathcal{P}_2\mathcal{P}_3$ &

all types &

$ {\rm I}_{1,1}^3$\\

\hline
%%%%%%%%%
$(0,0,0;1)$ & 

$\mathcal{P}_1\mathcal{P}_1^c\mathcal{P}_2\mathcal{P}_2^c\mathcal{P}_3\mathcal{P}_3^c$ &

all types &

$ (\ZZ/p\ZZ \times \mu_{p})^3$\\

\hline
%%%%%%%%%
$(1,0,0;1),\, (1,1,0;1),\, (0,1,0;1)$ & 

$\mathcal{P}_1\mathcal{P}_1^c\mathcal{P}_2\mathcal{P}_2^c\mathcal{P}_3$ &

all types &

$ (\ZZ/p\ZZ \times \mu_{p})^2 \times {\rm I}_{1,1}$\\

\hline
%%%%%%%%%
$(0,1,1;1),\, (0,1,0;1),\, (0,0,1;1)$ & 

$\mathcal{P}_1\mathcal{P}_1^c\mathcal{P}_2$ &

all types &

$ (\ZZ/p\ZZ \times \mu_{p})\times {\rm I}_{1,1}^2$\\

\hline
\end{tabular}
\end{table}

%%%%%%%%%%%%%%%%%%%%%%%%%%%%%

\begin{table}[h]
\caption{Elements of order $3$}
\begin{tabular}{|c|c|c|c|}

\hline
$\sigma$ & $p\cO_{K}$ & CM-types & $A[p] \times \bar{\FF}_p$\\
\hline
%%%%%%%%%
\multirow{2}{*}{$(0,0,0;s), \, (0,0,0;s^2)$} & 

\multirow{2}{*}{$\mathcal{P}\mathcal{P}^{c}$} &

(A), (B), (C) & 
${\rm I}_{3,2}$\\ \cline{3-4}

&
&
(D) &
$ (\ZZ/p\ZZ \times \mu_{p})^3$\\
\hline
%%%%%%%%%
\multirow{2}{*}{$(1,1,0;s), \, (0,1,1;s^2)$} & 

\multirow{2}{*}{$\mathcal{P}\mathcal{P}^{c}$} &

(A), (C), (D) & 
${\rm I}_{3,2}$\\ \cline{3-4}

&
&
(B) &
$ (\ZZ/p\ZZ \times \mu_{p})^3$\\
\hline
%%%%%%%%%
\multirow{2}{*}{$(1,0,1;s), \, (1,1,0;s^2)$} & 

\multirow{2}{*}{$\mathcal{P}\mathcal{P}^{c}$} &

(B), (C), (D) & 
${\rm I}_{3,2}$\\ \cline{3-4}

&
&
(A) &
$ (\ZZ/p\ZZ \times \mu_{p})^3$\\
\hline
\end{tabular}
\end{table}

%%%%%%%%%%%%%%%%%%%%%%%%%%%%%%%%%%

\begin{table}[h]
\caption{Elements of order $4$}
\begin{tabular}{|c|c|c|c|}
\hline
$\sigma$ & $p\cO_{K}$ & CM-types & $A[p] \times \bar{\FF}_p$\\
\hline
$(1,0,0;t),\, (0,1,0;t),\, (0,1,1;t), \, (1,0,1;t)$ & 

$\mathcal{P}_1\mathcal{P}_1^c\mathcal{P}_2$ &

all types &

$ (\ZZ/p\ZZ \times \mu_{p})\times {\rm I}_{2,1}$\\

\hline
$(1,0,0;ts),\, (0,0,0;ts),\, (0,1,1;ts), \, (1,1,0;ts)$ & 

$\mathcal{P}_1\mathcal{P}_1^c\mathcal{P}_2$ &

all types &

$ (\ZZ/p\ZZ \times \mu_{p})\times {\rm I}_{2,1}$\\
\hline
$(1,0,1;ts^2),\, (1,1,0;ts^2),\, (0,1,0;ts^2), \, (0,0,1;ts^2)$ & 

$\mathcal{P}_1\mathcal{P}_2$ &

all types &

$ {\rm I}_{2,1}\times {\rm I}_{1,1}$\\

\hline
\end{tabular}
\end{table}

%%%%%%%%%%%%%%%%%%%%%%%%%%%%%%%%%%

\begin{table}[h]
\caption{Elements of order $6$}
\begin{tabular}{|c|c|c|c|}
\hline
$\sigma$ & $p\cO_{K}$ & CM-types & $A[p] \times \bar{\FF}_p$\\
\hline
%%%%%%%%%
\multirow{2}{*}{$(0,0,1;s), (1,0,0;s^2)$} & 

\multirow{2}{*}{$\mathcal{P}$} &

(A), (B), (C) &

$ {\rm I}_{3,1}$\\ \cline{3-4}

&
&
(D) &
$ {\rm I}_{1,1}^3$\\
\hline
%%%%%%%%%
\multirow{2}{*}{$(1,1,1;s^2), (1,1,1;s)$} & 

\multirow{2}{*}{$\mathcal{P}$} &

(A), (B), (C) &

$ {\rm I}_{3,1}$\\ \cline{3-4}

&
&
(D) &
$ {\rm I}_{1,1}^3$\\
\hline
%%%%%%%%%
\multirow{2}{*}{$(1,0,0;s), (0,1,0;s^2)$} & 

\multirow{2}{*}{$\mathcal{P}$} &

(A), (B), (D) &

$ {\rm I}_{3,1}$\\ \cline{3-4}

&
&
(C) &
$ {\rm I}_{1,1}^3$\\
\hline

\end{tabular}
\end{table}
%%%%%%----------table---------%%%%%%%%%

As result we get the following corresponding
\newline
%\begin{table}[h]
\begin{tabular}{|c|c|c|c|}
\hline
ideal decomposition &$BT_{1}$ group scheme & $p$-rank & $a$-number\\
\hline
$P_{1} P_{1}^{c}P_{2} P_{2}^{c}P_{3} P_{3}^{c}$& $(\ZZ/p)^{3} \times (\mu_{p})^{3}$ & 3 & 0\\
\hline
\multirow{2}{*}{$P_{1} P_{1}^{c}P_{2} P_{2}^{c}$} &$(\ZZ/p)^{3} \times (\mu_{p})^3 $& 3&0\\\cline{2-4}
& $\ZZ/p \times \mu_{p} \times I_{1,1}^{2}$ & 1 &2\\
\hline
\multirow{2}{*}{$PP^{c}$}& $(\ZZ/p)^{3} \times (\mu_{p})^3 $ &3&0\\\cline{2-4}
&$I_{3,2}$ &0&2\\
\hline
$P_{1}P_{1}^{c}P_{2}P_{2}^{c}P_{3}$ &$(\ZZ/p\times \mu_{p})^{2}\times I_{1,1}$ &2&1\\
\hline
\multirow{2}{*}{$P_{1}P_{1}^{c}P_{2}$}& $(\ZZ/p \times \mu_{p})^{2} \times I_{1,1}$ & 2&1\\\cline{2-4}
& $I_{1,1}^{3}$ &0&3\\
\hline
\multirow{2}{*}{$P_{1}P_{1}^{c}P_{2}P_{3}$} &$(\ZZ/p)^{3} \times (\mu_{p})^{3}$& 3&0\\\cline{2-4}
& $\ZZ/p \times \mu_{p} \times I_{1,1}^{2}$ &1&2\\
\hline
$P_{1}P_{2}P_{3}$ & $I_{1,1}^{3}$ &0&3\\
\hline
$P_{1}P_{2}$ & $I_{2,1} \times I_{1,1}$ & 0&2\\
\hline
\multirow{2}{*}{$P$} & $I_{3,1}$ &0 &1\\\cline{2-4}
& $I_{1,1}^{3}$& 0&3\\
\hline
\end{tabular}
%\end{table}
\end{subsubsection}

Similarly we get the following results
\begin{subsubsection}{}${\rm Gal}(\tilde{K}/K) \cong \ZZ_{6}$

%If ${\rm Gal}(\tilde{K}/K) \cong \ZZ_{6}$ then the following relations hold
\begin{tabular}{|c|c|c|c|}
\hline
ideal decomposition &$BT_{1}$ group schemes & $p$-rank & $a$-number\\
\hline
$P_{1} P_{1}^{c}P_{2} P_{2}^{c}P_{3} P_{3}^{c}$& $(\ZZ/p)^{3} \times (\mu_{p})^{3}$ & 3 & 0\\
\hline
\multirow{2}{*}{$PP^{c}$}& $(\ZZ/p)^{3} \times (\mu_{p})^3 $ &3&0\\\cline{2-4}
&$I_{3,2}$ &0&2\\
\hline
$P$ & $I_{3,1}$ &0 &1\\
\hline
\end{tabular}

\end{subsubsection}
%%%%%%%%%%%%%%%
\begin{subsubsection}{}{${\rm Gal}(\tilde{K}/\QQ) \cong \ZZ_{2} \times S_{3}$}

%If ${\rm Gal}(\tilde{K}/\QQ) \cong \ZZ_{2} \times S_{3}$ then the following relations hold

\begin{tabular}{|c|c|c|c|}
\hline
ideal decomposition &$BT_{1}$ group schemes & $p$-rank & $a$-number\\
\hline
$P_{1} P_{1}^{c}P_{2} P_{2}^{c}P_{3} P_{3}^{c}$& $(\ZZ/p)^{3} \times (\mu_{p})^{3}$ & 3 & 0\\
\hline
\multirow{2}{*}{$P_{1} P_{1}^{c}P_{2} P_{2}^{c}$} &$(\ZZ/p)^{3} \times (\mu_{p})^3 $& 3&0\\\cline{2-4}
& $\ZZ/p \times \mu_{p} \times I_{1,1}^{2}$ & 1 &2\\
\hline
\multirow{2}{*}{$PP^{c}$}& $(\ZZ/p)^{3} \times (\mu_{p})^3 $ &3&0\\\cline{2-4}
 &$I_{3,2}$ &0&2\\
\hline
\multirow{2}{*}{$P_{1}P_{1}^{c}P_{2}$}& $(\ZZ/p \times \mu_{p})^{2} \times I_{1,1}$ & 2&1\\\cline{2-4}
& $I_{1,1}^{3}$ &0&3\\
\hline
$P_{1}P_{2}P_{3}$ & $I_{1,1}^{3}$ &0&3\\
\hline
$P_{1}P_{2}$ & $I_{2,1} \times I_{1,1}$ & 0&2\\
\hline
$P$ & $I_{3,1}$ &0 &1\\
\hline
\end{tabular}
\end{subsubsection}
%%%%%%%%%%%%%%%%

\begin{subsubsection}{}{${\rm Gal}(\tilde{K}/\QQ) \cong (\ZZ_{2})^{3} \rtimes \ZZ_{3}$}
%If ${\rm Gal}(\tilde{K}/\QQ) \cong (\ZZ_{2})^{3} \rtimes \ZZ_{3}$ then the following relations hold

\begin{tabular}{|c|c|c|c|}
\hline
ideal decomposition &$BT_{1}$ group schemes & $p$-rank & $a$-number\\
\hline
$P_{1} P_{1}^{c}P_{2} P_{2}^{c}P_{3} P_{3}^{c}$& $(\ZZ/p)^{3} \times (\mu_{p})^{3}$ & 3 & 0\\
\hline
\multirow{2}{*}{$PP^{c}$}& $(\ZZ/p)^{3} \times (\mu_{p})^3 $ &3&0\\\cline{2-4}
 &$I_{3,2}$&0&2\\
\hline
$P_{1}P_{1}^{c}P_{2}P_{2}^{c}P_{3}$ &$(\ZZ/p\times \mu_{p})^{2}\times I_{1,1}$ &2&1\\
\hline
\multirow{2}{*}{$P_{1}P_{1}^{c}P_{2}P_{3}$} &$(\ZZ/p)^{3} \times (\mu_{p})^{3}$& 3&0\\\cline{2-4}
& $\ZZ/p \times \mu_{p} \times I_{1,1}^{2}$ &1&2\\
\hline
$P_{1}P_{2}P_{3}$ & $I_{1,1}^{3}$ &0&3\\
\hline
\multirow{2}{*}{$P$} & $I_{3,1}$ &0 &1\\\cline{2-4}
& $I_{1,1}^{3}$& 0&3\\
\hline
\end{tabular}

\end{subsubsection}
 
\end{subsection}
\end{section}

\begin{section}{speculation about new parameter, $b$-number, and densities}
Let $A$ be an absolutely simple ordinary Abelian variety over  $\bar{\FF}_p$, where $p$ is a prime number. Then an algebra
$K:={\rm End}(A)\otimes \QQ$ is a CM field with $[K:\QQ]=2\, {\rm dim}(A)$. Let $K_{0}$ be a maximal totally real subfield of $K$ and $\cO_{K_0}$ be its full ring of integers. Then we say that $m$ is {\bf $b-$number} of  $A$, denoted by {\bf $b(A)$}, if
$$
p\cO_{K_0}=\mathcal{P}_1^{i_1} \ldots \mathcal{P}_m^{i_m}.
$$ 
If $A$ is isogenous over $\bar{\FF}_p$ to a product of absolutely simple Abelian varieties $A_{1}\times \ldots \times A_t$ then we define 
\begin{displaymath}
b(A)=\sum_{i=1}^{t} b(A_i).
\end{displaymath}

In the same way we define {\bf $b(\FF_q)-$number}, namely
let $A$ be a simple ordinary Abelian variety over $\FF_q$, with $p={\rm char}( \FF_q)$. Then an algebra
$K:={\rm End}_{\FF_q}(A)\otimes \QQ$ is a CM field with $[K:\QQ]=2\, {\rm dim}(A)$. Let $K_{0}$ be a maximal totally real subfield of $K$ and $\cO_{K_0}$ be its full ring of integers. Then we say that $m$ is {\bf $b(\FF_q)-$number} of  $A$ if
$$
p\cO_{K_0}=\mathcal{P}_1^{i_1} \ldots \mathcal{P}_m^{i_m}.
$$ 
If $A$ is isogenous over $\FF_q$ to a product of simple Abelian varieties $A_{1}\times \ldots \times A_t$ then we define 
\begin{displaymath}
b(\FF_q)(A)=\sum_{i=1}^{t} b(\FF_q)(A_i).
\end{displaymath}

\begin{remark}
The following  properties of $b-$number hold,
\begin{itemize}
\item{if Abelian varieties $A$ and $B$ are $\bar{\FF}_p$ isogenous then $b(A)=b(B)$,}
\item{if Abelian varieties $A$ and $B$ are $\FF_q$ isogenous then $b(\FF_q)(A)=b(\FF_q)(B)$, }
\item{$1 \le b(A) \le {\rm dim}(A)$, the same holds for $b(\FF_q)(A).$}
\end{itemize}
\end{remark}

We can define a subset of $\bar{\FF}_p$--isogeny classes  
\[
B_{i,g}=\{[A] \in\mathcal{A}_{g,1}^{ord} \otimes \FF_p/\sim \quad | \quad b(A)=i \}
\]
of set   all $\bar{\FF}_p$--isogeny classes of  ordinary principally polarized Abelian varieties $\mathcal{A}_{g,1}^{ord} \otimes \FF_p/\sim$.

Similarly, One can define a subset  of $\mathcal{A}_{g,1}^{ord}(\FF_q)/\sim$, the set of the isogeny classes of ordinary principally polarized Abelian varieties over $\FF_q$,   as 
\[
B_{i,g}(\FF_q)=\{[A] \in\mathcal{A}_{g,1}^{ord}(\FF_q)/\sim \quad | \quad b(\FF_q)(A)=i \}.
\]

Unfortunately, the geometric structure of $B_{i,g}$ and $B_{i,g}(\FF_q)$ is not clear and might be that it does not exits at all. However, we can introduce a notion of density 
\[
{\rm D}(B_{i,g}(\FF_q))=\frac{\# B_{i,g}(\FF_q)}{\#\left( \mathcal{A}_{g,1}^{ord}(\FF_q)/\sim\right)}
\]
of the subset $B_{i,g}(\FF_q)$ of the set of all isogeny classes of 
ordinary principally polarized Abelian varieties over $\FF_q$. Note that this definition is well defined since the set $\mathcal{A}_{g,1}^{ord}(\FF_q)/\sim$ is finite.

Now we can define a density

\[
{\rm D}(B_{i,g})=\lim_{n \to \infty}\frac{\# B_{i,g}(\FF_{p^n})}{\#\left( \mathcal{A}_{g,1}^{ord}(\FF_{p^n})/\sim\right)}
\]
of the subset $B_{i,g}(\FF_q)$ of the set of all isogeny classes of 
ordinary principally polarized Abelian varieties over $\bar{\FF}_p$ under condition that the limit
exits.

Here we give some informal speculation about connection of this density with number theory.
According to the Honda-Tate theory it seems that the density ${\rm D}(B_{m,g})$ is closely related to the number
\[
\lim_{\delta \to \infty}
\frac
{\#\{K|\, \mbox{totally real,}\, [K:\QQ]=g,\, p\cO_{K}=\mathcal{P}_1^{i_1} \ldots \mathcal{P}_m^{i_m},\,
\mbox{with discriminate}\, \delta(K) \le \delta\}}
{\#\{K|\, \mbox{totally real,}\, [K:\QQ]=g,\, \delta(K) \le \delta\}}
\]
which looks like the Chebotarive density theorem but for fields. In other words, it tells about a density of real fields with prescribed decomposition type of a fixed prime number $p$ and fixed extension degree among all real fields with the same extension degree.

Further we produce an example  which is based on  the following lemma.
\begin{lemma}
A monic quadratic polynomial in $\ZZ[X]$
\[
X^4+aX^3+bX^2+qaX+q^2
\]
is a Weil polynomial if and only if 
$|a| \le 4\sqrt{q}$ and  $2|a|\sqrt{q}-2\sqrt{q} \le b \le a^2/4 +2q.$ Moreover it corresponds to an ordinary Abelian surface if and only if $b \not \equiv 0\, {\rm mod}\, p.$
\end{lemma}

\begin{example}
We adopt the lemma above  for  polynomials of the form $T^2+aT+(b-2q)$ and $p>2$. As result  the following hold
\begin{itemize}
\item{$p$ splits completely if 
$\left(\frac{\Delta}{p}\right)=1$,}
\item{ $p$ is inert if $\left(\frac{\Delta}{p}\right)=-1$,}
\item{ $p$ is ramfied if $\left(\frac{\Delta}{p}\right)=0$,}
\end{itemize}
where $\Delta=a^2-4(b-2q)$.

 So for $p=3$ and $p=5$ we have

\begin{tabular}{|c|c|c|}
\hline
$\FF_q$ & ${\rm D}(B_{1,1})$ &${\rm D}(B_{1,2})$\\
\hline
$\FF_3$ & $0.3388888889$ &$0.6611111111$\\
\hline
$\FF_{3^2}$ & $0.3383233533$ &$0.6616766467$\\
\hline
$\FF_{3^3}$ & $0.3342293907$ &$0.6657706093$\\
\hline
$\FF_{3^5}$ & $0.3334656710$ &$0.6665343290$\\
\hline
$\FF_{3^8}$ & $0.3333404746$ &$0.6666595254$\\
\hline
\end{tabular}

\phantom{space}
%\vspace{0.5 cm}
\begin{tabular}{|c|c|c|}
\hline
$\FF_q$ & ${\rm D}(B_{1,1})$ &${\rm D}(B_{1,2})$\\
\hline
$\FF_5$ & $0.4070080863$ &$0.5929919137$\\
\hline
$\FF_{5^2}$ & $0.4019172317$ &$0.5980827683$\\
\hline
$\FF_{5^3}$ & $0.4001517555$ &$0.5998482445$\\
\hline
$\FF_{5^5}$ & $0.4000055751$ &$0.5999944249$\\
\hline
\end{tabular}

\end{example}

\end{section}

%------------------------------------------------------------------------------------

\bibliographystyle{plain}
\def\cprime{$'$}

\end{document}